\documentclass[10pt]{article}\usepackage{pb-diagram}

\usepackage{amssymb,latexsym}
\usepackage{amsfonts}
\usepackage{amsmath,amsthm,graphicx}

\newcommand{\bd}{\begin{document}}
\newcommand{\ed}{\end{document}}
\newcommand{\bc}{\begin{center}}
\newcommand{\ec}{\end{center}}
\newcommand{\vs}{\vspace}
\newcommand{\hs}{\hspace}
\newcommand{\bq}{\begin{quote}}
\newcommand{\eq}{\end{quote}}
\newcommand{\mb}{\makebox}
\newcommand{\lt}{\left}
\newcommand{\rt}{\right}
\newcommand{\beqa}{\begin{eqnarray*}}\large
\newcommand{\eeqa}{\end{eqnarray*}}
\newcommand{\beqn}{\begin{eqnarray}}
\newcommand{\eeqn}{\end{eqnarray}}
\newcommand{\bbibl}{\begin{thebibliography}{99}}
\newcommand{\ebibl}{\end{thebibliography}}
\newcommand{\ti}{\times}
\newcommand{\bit}{\begin{itemize}}
\newcommand{\eit}{\end{itemize}}
\newcommand{\ben}{\begin{enumerate}}
\newcommand{\een}{\end{enumerate}}
\newcommand{\lb}{\label}
\newcommand{\hf}{\hspace*{\fill}}
\newcommand{\vf}{\vspace*{\fill}}
\newcommand{\beq}{\begin{equation}}
\newcommand{\eeq}{\end{equation}}
\newcommand{\ba}{\begin{array}}
\newcommand{\ea}{\end{array}}
\newcommand{\del}{\partial}
\newcommand{\bm}[1]{\mb{\boldmath ${#1}$}}
\newcommand{\ot}{\otimes}
\newcommand{\nn}{\nonumber}
\newcommand{\R}{\mb{$I\!\!R$}}
\newcommand{\C}{{\cal C}}
\newcommand{\M}{{\cal M}}
\newcommand{\E}{{\cal E}}
\newcommand{\N}{{\cal N}}
\newcommand{\B}{{\cal B}}
\newcommand{\Y}{{\cal Y}}
\newcommand{\F}{{\cal F}}
\newcommand{\Rc}{{\cal R}}
\newcommand{\A}{{\cal A}}
\renewcommand{\P}{{\cal P}}
\renewcommand{\S}{{\cal S}}
\newcommand{\es}{\emptyset}
\newcommand{\ci}{\subseteq}
\newcommand{\cs}{\supseteq}
\renewcommand{\u}{\cup}
\renewcommand{\i}{\cap}
\newcommand{\bu}{\bigcup}
\newcommand{\bi}{\bigcap}
\newcommand{\la}{\leftarrow}
\newcommand{\ra}{\rightarrow}
\newcommand{\Ra}{\Rightarrow}
\newcommand{\Lra}{\Leftrightarrow}
\newcommand{\lgra}{\longrightarrow}
\newcommand{\Lgra}{\Longrightarrow}
\newcommand{\lglra}{\longleftrightarrow}
\newcommand{\Lglra}{\Longleftrightarrow}
\renewcommand{\a}{\alpha}
\renewcommand{\b}{\beta}
\newcommand{\g}{\gamma}
\newcommand{\G}{\Gamma}
\renewcommand{\d}{\delta}
\newcommand{\D}{\Delta}
\newcommand{\e}{\varepsilon}
\newcommand{\eps}{\epsilon}
\newcommand{\h}{\eta}
\renewcommand{\l}{\lambda}
\newcommand{\m}{\mu}
\newcommand{\n}{\nu}
\newcommand{\p}{\pi}
\newcommand{\s}{\sigma}
\newcommand{\Si}{\Sigma}
\newcommand{\ta}{\tau}
\newcommand{\ph}{\phi}
\newcommand{\Ph}{\Phi}
\renewcommand{\c}{\chi}
\newcommand{\om}{\omega}
\newcommand{\Om}{\Omega}
\newcommand{\tri}{\triangle}
\newcommand{\rec}[1]{\frac{1}{#1}}
\newcommand{\f}{\frac}
\newcommand{\sm}[2]{\sum_{#1}^{#2}}
\newcommand{\ld}{\ldots}
\newcommand{\ov}{\overline}
\newcommand{\ol}[1]{$\bar{\mb{#1}}$}
\newcommand{\un}{\underline}
\newcommand{\iy}{\infty}

\newcommand{\wt}{\widetilde}
\newcommand{\ds}{\displaystyle}
\newcommand{\bdm}{\begin{displaymath}}
\newcommand{\edm}{\end{displaymath}}

\newcommand{\nin}{\not\in}
\newcommand{\bt}{\begin{tabular}}
\newcommand{\et}{\end{tabular}}
\newcommand{\alter}[2]{\lt\{ \ba {ll}#1 \\ #2 \ea \rt.}
\newcommand{\alt}[4]{\lt\{ \ba{ll}#1 & \mb{if \,\,}#2 \\ #3 & \mb{if
               \,\,}#4 \ea \rt.}
\newcommand{\altn}[4]{\lt\{ \ba{rl}#1 & \mb{if \,\,}#2 \\ #3 & \mb{if
               \,\,}#4 \ea \rt.}
\newcommand{\alto}[6]{ \lt\{ \ba{ll}#1 & \mb{if \,\,}#2 \\ #3 & \mb{if
               \,\,} #4 \\ #5 & \mb{if \,\,}#6 \ea \rt.}
\newcommand{\altero}[5]{\mb{$\lt\{ \ba {ll}#1 & \mb{if \,\,}#2 \\ #3 &
               \mb{if \,\,} #4 \\ #5 & \mb{otherwise} \ea \rt.$}}

\newcounter{cnt1}
\newcounter{cnt2}
\newcounter{cnt3}
\newcommand{\blr}{\begin{list}{$($\roman{cnt1}$)$} {\usecounter{cnt1}
        \setlength{\topsep}{0pt} \setlength{\itemsep}{0pt}}}
\newcommand{\bla}{\begin{list}{$($\alph{cnt2}$)$} {\usecounter{cnt2}
        \setlength{\topsep}{0pt} \setlength{\itemsep}{0pt}}}
\newcommand{\bln}{\begin{list}{$($\arabic{cnt3}$)$} {\usecounter{cnt3}
                \setlength{\topsep}{0pt} \setlength{\itemsep}{0pt}}}
\newcommand{\el}{\end{list}}
\newcommand{\no}{\noindent}
\newtheorem{Thm}{Theorem}[section]
\newtheorem{Lem}[Thm]{Lemma}
\newtheorem{Prop}[Thm]{Proposition}
\newtheorem{Def}[Thm]{Definition}
\newtheorem{Exm}[Thm]{Example}
\newtheorem{Rem}[Thm]{Remark}
\newtheorem{Cor}[Thm]{Corollory}
\renewcommand{\baselinestretch}{1}
\newcommand{\ilim}{\mathop{\varprojlim}\limits}
\newcommand{\dlim}{\mathop{\varinjlim}\limits}

\begin{document}
\title{On a theorem of Eilenberg in simplicial Bredon-Illman cohomology with local coefficients}
\author{Goutam Mukherjee and Debasis Sen}
\date{}
\maketitle{}
\noindent

\begin{abstract}{ We prove simplicial version of a classical theorem of Eilenberg in the equivariant context and give an alternative description of  the simplicial version of Bredon-Illman cohomology with local coefficients, as introduced in\cite{ms}, to derive a spectral sequence.}
\end{abstract}
\footnote {\bf The second author would like to thank CSIR for its support.}
\section{Introduction} For spaces with group actions the analogue of cohomology with local coefficients is the Bredon-Illman cohomology with local coefficients, as introduced in \cite{mm}. This is based on the notion of fundamental groupoid of a space equipped with a group action.  A classical theorem of Eilenberg states that cohomology with local coefficients of a space can be described by the cohomology of an invariant subcomplex of the cochain complex of its universal cover, where the universal cover is equipped with the action of the fundamental group of the base space. An equivariant analogue of this result was proved in \cite{mm}. Recently, in \cite{ms}, we introduced  equivariant simplicial cohomology with local coefficients, which is the simplicial version of Bredon-Illman cohomology with local coefficients and proved a classification theorem. The corresponding non-equivariant result was proved in \cite{hir}, \cite{gj}, \cite{bfg}. In this paper we derive Eilenberg's theorem for equivariant simplicial cohomology with local coefficients. This is based on the notion of universal covering complexes of one vertex Kan complexes \cite{gugg}.  In equivariant context, the role of the universal cover is played by a contravariant functor from the category of canonical orbits to the category of one vertex Kan complexes. Finally,
we give an alternative description of equivariant simplicial cohomology with local coefficients via the notion of cohomology of a small category following \cite{mp} and use it to derive a spectral sequence.

The paper is organized as follows. In Section $2$, we recall some standard results and fix notations. The notion of equivariant local coefficients of a simplicial set equipped with a simplicial group action is based on fundamental groupoid. In Section 3, we recall these concepts and the definition of simplicial version of Bredon-Illman cohomology with local coeffic

\section{Preliminaries on $G$-simplicial sets}

In this section we set up our notations and recall some standard facts \cite{may}.

Let $\Delta$ be the category whose objects are ordered sets
$$[n]=\{0 < 1< \cdots <n\},~n\geq  0,$$ and morphisms are non-decreasing maps $f : [n] \lgra [m].$ There are some distinguished morphisms $d^i:[n-1]\lgra [n], 0\leq i \leq n$, called cofaces and $s^i : [n+1] \lgra [n],~ 0\leq i \leq n$, called codegeneracies, defined as follows:
$$d^i(j) = j,~j<i~~\mbox{and}~~ d^i(j) = j+1,~ j\geq i, ~~(n>0,~  0\leq i \leq n);$$
$$s^i(j) = j,~j \leq i,~~\mbox{and}~~s^i(j) = j-1,~j> i,~~(n \geq 0, ~  0\leq i \leq n).$$
These maps satisfy the standard cosimplicial relations.

A simplicial object $X$ in a category $\mathcal{C}$ is a contravariant functor $X: \Delta \lgra \mathcal{C}.$ Equivalently, a simplicial object is a sequence $\{X_n\}_{n\geq 0}$ of objects of $\mathcal{C}$, together with $\mathcal{C}$-morphisms $ \partial_i : X_n\lgra X_{n-1}$ and $s_i:X_n\lgra X_{n+1},$ $0\leq i \leq n,$
verifying the following simplicial identities:
$$\partial_i \partial_j = \partial_{j-1} \partial_i,~ ~  \partial_i s_j = s_{j-1} \partial_i,~ \mbox{if}~~i<j,$$
$$~~\partial_j s_j =id = \partial_{j+1}s_j,$$
$$\partial_i s_j = s_j \partial_{i-1}, ~~i >j+1;~~ s_is_j = s_{j+1} s_i,~~i\leq j.$$
A simplicial map $f:X\lgra Y$ between two simplicial objects in a category $\mathcal{C}$, is a collection of $\mathcal{C}$-morphisms $f_n : X_n\lgra Y_n,$ $n\geq 0,$  commuting with $\partial_i$ and $s_i$.

In particular, a simplicial set is a simplicial object in the category of sets. Throughout $\mathcal{S}$ will denote the category of simplicial sets and simplicial maps.

For any $n$-simplex $x\in X_n$, in a simplicial set $X$, we shall use the notation $\partial_{(i_1, i_2, \cdots, i_r)}x$ to denote the simplex $\partial_{i_1}\partial_{i_2}\cdots \partial_{i_r}x$ obtained by applying the successive face maps $\partial_{i_{r-k}}$ on $x$, where $0\leq i_{r-k}\leq n-k,~0\leq k \leq r-1.$

Recall that the simplicial set $\Delta[n]$, $n\geq 0$, is defined as follows. The set of $q$-simplices is
$$\Delta [n]_q = \{(a_0,a_1,\cdots,a_q) ;~~\mbox{where}~ a_i\in \mathbb{Z},~ 0\leq a_0\leq a_1\leq \cdots \leq a_q\leq n\}.$$
The face and degeneracy maps are defined by $$\partial_i(a_0,\cdots ,a_q)=(a_0,\cdots ,a_{i-1},a_{i+1},\cdots ,a_q),~ s_j(a_0,\cdots ,a_q)=(a_0,\cdots ,a_i,a_i,\cdots ,a_q).$$
Alternatively, the set of $k$-simplices  can be viewed as the contravariant functor
$$\Delta [n]([k]) = ~\mbox{Hom}_{\Delta}([k], [n]),$$
the set of $\Delta$-morphisms from $[k]$ to $[n]$. The only non-degenerate $n$-simplex is $id : [n] \lgra [n]$ and is denoted by $\Delta_n$. In the earlier notation, it is simply, $\Delta_n = (0, 1, \cdots, n).$

It is well known that if $X$ is a  simplicial set, then for any $n$-simplex $x\in X_n$ there is a unique simplicial map $\overline{x} : \Delta[n]\lgra X$ with $\overline{x}(\Delta_n) =x$. Often by an $n$-simplex in a simplicial set $X$ we shall mean either an element $x \in X_n$ or the corresponding simplicial map $\overline{x}.$

We have simplicial maps $\delta_i:\Delta[n-1]\rightarrow \Delta[n]$ and $\sigma_i:\Delta[n+1]\rightarrow \Delta[n]$ for $0\leq i \leq n$ defined by $\delta_i(\Delta_{n-1})=\partial_i(\Delta_n)$ and $\sigma_i(\Delta_{n+1})=s_i(\Delta_n)$. The boundary subcomplex $\partial \Delta[n]$ of $\Delta[n]$ is defined as the smallest subcomplex of $\Delta[n]$ containing the faces $\partial_i\Delta_n,~~ i=0,1,...,n $.

\begin{Def}
Let $G$ be a discrete group.  A $G$-simplicial set is a simplicial object in the category of $G$-sets. More precisely, a $G$-simplicial set is a simplicial set $\{X_n ; \partial_i, s_i, 0\leq i \leq n\}_{n\geq 0}$ such that each $X_n$ is a $G$-set and the face maps $\partial_i: X_n\longrightarrow X_{n-1}$ and the degeneracy maps $s_i: X_n\longrightarrow X_{n+1}$ commute with the $G$-action.
A map between G-simplicial sets is a simplicial map which commutes with the G-action.
\end{Def}

\begin{Def}
A $G$-simplicial set $X$ is called $G$-connected if each fixed point simplicial set $X^H$, $H\subseteq G$, is connected.
\end{Def}

\begin{Def}
Two $G$-maps $f,g:K\rightarrow L$ between two $G$-simplicial sets are $G$-homotopic if there exists a $G$-map $F:K\times \Delta[1]\rightarrow L$ such that
$$F\circ (id\times \delta_1) =f,~~F\circ (id\times \delta_0) =g.$$
The map $F$ is called a $G$-homotopy from $f$ to $g$ and we write $F:f\simeq_G g.$
If $i:K^{\prime}\subseteq K$ is an inclusion of subcomplex and $f,~g$ agree on $K^{\prime}$ then we say that $f$ is $G$-homotopic to $g$ relative to $K^{\prime}$ if there exists a $G$ homotopy $F:f\simeq_G g$ such that $F\circ (i\times id) =\a \circ pr_1,$
where $\alpha=f|_K^{\prime}=g|_K^{\prime}$ and $pr_1 :K^{\prime}\times \Delta[1] \lgra K^{\prime}$ is the projection onto the first factor. In this case we write $F:f\simeq_G g(rel~ K^{\prime}).$
\end{Def}

\begin{Def}
A $G$-simplicial set is a $G$-Kan complex if for every subgroup $H\subseteq G$ the fixed point simplicial set $X^H$ is a Kan complex.
\end{Def}
\begin{Rem}\lb{equi}
Recall (\cite{ag}, \cite{fg}) that the category $G\mathcal{S}$ of $G$-simplicial sets and $G$-simplicial maps between $G$-simplicial sets has a closed model structure \cite{qui}, where the fibrant objects are the $G$-Kan complexes and cofibrant objects are the $G$-simplicial sets. From this it follows that $G$-homotopy on the set of $G$-simplicial maps $K\rightarrow L$ is an equivalence relation, for every $G$-simplicial set $K$ and $G$-Kan complex $L$. More generally, relative $G$-homotopy is an equivalence relation if the target is a $G$-Kan complex.
\end{Rem}

We consider $G/H \times \Delta[n]$ as a simplicial set where $(G/H \times \Delta[n])_q=G/H \times (\Delta[n])_q$ with face and degeneracy maps as $id\times \partial_i$ and $id \times s_i$. Note that the group $G$ acts on $G/H$  by left translation. With this $G$-action  on the first factor and trivial action on the second factor $G/H \times \Delta[n]$ is a $G$-simplicial set.

A $G$-simplicial map $\sigma:G/H \times \Delta[n]\rightarrow X$ is called an equivariant $n$-simplex of type $H$ in $X$.
\begin{Rem}\lb{correspondence}
We remark that for a $G$-simplicial set $X,$ the set of equivariant $n$-simplices in $X$ is in bijective correspondence with n-simplices of $X^H$. For an equivariant $n$-simplex $\sigma$, the corresponding $n$-simplex is $\sigma^{\prime}=\sigma (eH, \Delta_n).$ The simplicial map $\Delta[n]\lgra X^H,~~\Delta_n \mapsto \sigma^{\prime}$
will be denoted by $\overline{\sigma}.$
\end{Rem}

We shall call $\sigma$  degenerate or non-degenerate according as the $n$-simplex $\sigma^{\prime} \in X^H_n$ is degenerate or non-degenerate.

Recall that the category of canonical orbits, denoted by $O_G,$ is a category whose  objects are cosets $G/H$, as $H$ runs over the all subgroups of $G$. A morphism from $G/H$ to $G/ K$ is a $G$-map. Recall that such a morphism determines and is determined by a subconjugacy relation $g^{-1}Hg\subseteq K$ and is given by $\hat{g}(eH)=gK$. We denote this morphism by $\hat{g}$ \cite{br}.
\begin{Def}
A contravariant functor from $O_G$ to $\mathcal{S}$ is called an $O_G$-simplicial set. A map between $O_G$-simplicial sets is a natural transformation of functors.
\end{Def}
We shall denote the category of $O_G$-simplicial sets by $O_G\mathcal{S}.$

The notion of $O_G$-groups or $O_G$-abelian groups has the obvious meaning replacing $\mathcal{S}$ by $\mathcal{G}rp$ or $\mathcal{A}b.$

For a $G$-simplicial set $X$, with a $G$-fixed $0$-simplex $v$, we have an $O_G$-group $\underline{\pi}X$ defined as follows. For any subgroup $H$ of $G$,
$$\underline{\pi}X(G/H) := \pi_1(X^H,v)$$ and for a morphism $\hat{g}:G/H\lgra G/K,~~g^{-1}Hg\subseteq K$, $\underline{\pi}X(\hat{g})$ is the homomorphism in fundamental groups induced by the simplicial map $g:X^K\lgra X^H.$
\begin{Def}
An $O_G$-group $\underline{\pi}$ is said to act on an $O_G$-simplicial set (group or abelian group) $\underline{X}$ if for every subgroup $H\subseteq G$, $\underline{\pi}(G/H)$ acts on $\underline{X}(G/H)$ and this action is natural with respect to maps of $O_G.$ Thus if
$$\phi (G/H):\underline{\pi}(G/H)\times \underline{X}(G/H)\lgra \underline{X}(G/H)$$ denotes the action of $\underline{\pi}(G/H)$ on $\underline{X}(G/H)$ then for each subconjugacy relation\\ $g^{-1}Hg\subseteq K,$
$$\phi (G/H)\circ (\underline{\pi}(\hat{g})\times\underline{X}(\widehat{g}))= \underline{X}(\hat{g})\circ \phi(G/K).$$
\end{Def}

\section{Simplicial Bredon-Illman Cohomology with local coefficients} In this section we recall \cite{ms} the notion of fundamental groupoid of a $G$-simplicial set $X$, the notion of equivariant local coefficients on $X$ and
the definition of simplicial Bredon-Illman cohomology with local coefficients.

We begin with the notion of fundamental groupoid. Recall \cite{gj} that the fundamental groupoid $\pi X$ of a Kan complex $X$ is a category having as objects all $0$-simplexes of $X$ and a morphism $ x\longrightarrow y$ in $\pi X$ is a homotopy class of $1$-simplices $\omega : \Delta [1] \longrightarrow X$ rel $\partial \Delta [1]$ such that $\omega \circ \delta_0 = \overline{y}$, $\omega \circ \delta_1 = \overline{x}$. If $\omega_2$ represents an arrow from $x$ to $y$ and $\omega_0$ represents an arrow from $y$ to $z$,  then their composite $[\omega_0]\circ [\omega_2]$ is represented by $\Omega \circ \delta_1$, where the simplicial map $\Omega:\Delta[2]\lgra X$ corresponds to a $2$-simplex, which is determined by the compatible pair $(\omega^{\prime}_0,~~, \omega^{\prime}_2)$. For a simplicial set $X$ the notion of fundamental groupoid is defined via the geometric realization and the total singular functor.

Suppose $x_H$ and $y_K$ are equivariant $0$-simplices of type $H$ and $K$, respectively, and $\hat{g}:G/H \rightarrow G/K$ is a morphism in $O_G$, given by a subconjugacy relation $g^{-1}Hg\subseteq K$, $g\in G,$ so that $\hat{g}(eH)=gK$. Moreover suppose that we have an equivariant $1$-simplex $\phi:G/H \times \Delta[1]\rightarrow X$  of type $H$ such that
$$\phi \circ (id\times \delta_1) = x_H,~~\phi \circ (id\times \delta_0)=y_K\circ (\hat{g}\times id).$$
Then, in particular, $\phi^{\prime}$ is a $1$-simplex in $X^H$  such that $\partial_1\phi^{\prime} = x_H^{\prime}$ and $\partial_0\phi^{\prime} = gy_K^{\prime}$, notations are as in Remark \ref{correspondence}. Observe that the $0$-simplex $gy_K^{\prime}$ in $X^H$ corresponds to the composition
$$G/H \times \Delta [0]\stackrel{\hat{g}\times id}{\rightarrow} G/K \times \Delta [0] \stackrel{y_K}{\longrightarrow} X$$ and $\phi $ is a $G$-homotopy $x_H \simeq_G y_K \circ (\hat{g}\times id)$.

\begin{Def}
Let $X$ be a $G$-Kan complex. The fundamental groupoid $\Pi X$ is a category with objects equivariant $0$-simplices
$$x_{H}:G/ H \times \Delta[0] \rightarrow X$$
of type $H$, as $H$ varies over all subgroups of $G$. Given two objects $x_H$ and  $y_K$ in $\Pi X$, a morphism from $x_H \longrightarrow y_K$ is defined as follows. Consider the set of all pairs $(\hat{g},\phi)$ where $\hat{g}:G/H \rightarrow G/K$ is a morphism in $O_G$, given by a subconjugacy relation $g^{-1}Hg\subseteq K$, $g\in G$ so that $\hat{g}(eH)=gK$ and $\phi:G/H \times \Delta[1]\rightarrow X$ is an equivariant $1$-simplex such that
$$\phi \circ (id\times \delta_1) = x_H,~~\phi \circ (id\times \delta_0)=y_K\circ (\hat{g}\times id).$$

The set of morphisms in $\Pi X$ from $x_H$ to $y_K$ is a quotient of the set of pairs mentioned above by an equivalence relation $`\sim `,$ where $(\hat{g}_{1},\phi_{1})\sim(\hat{g}_{2},\phi_{2})$ if and only if  $g_1=g_2=g$ (say) and there exists a $G$-homotopy
$\Phi : G/H \times \Delta [1] \times \Delta [1] \longrightarrow X$ of $G$-homotopies such that $\Phi : \phi_1 \simeq_G \phi_2$ (rel $G/H \times \partial \Delta [1]$). Since $X$ is a $G$-Kan complex, by Remark \ref{equi}, $\sim$ is an equivalence relation. We denote the equivalence class of $(\hat{g},\phi)$ by $[\hat{g},\phi]$. The set of equivalence classes is the set of morphisms in $\Pi X$ from $x_H$ to $y_K$.

The composition of morphisms in $\Pi X$ is defined as follows. Given two morphisms
$$
\begin{diagram}
\node{x_{H}} \arrow{e,t}{[\hat{g}_{1},\phi_{1}]} \node{y_{K}} \arrow{e,t}{[\hat{g}_{2},\phi_{2}]} \node{z_{L}}
\end{diagram}
$$
their composition $[\hat{g}_2, \phi_2]\circ [\hat{g}_1, \phi_1]$ is $[\widehat{g_1g_2}, \psi]: x_H \longrightarrow z_L$, where the first factor is the composition
\[
\begin{diagram}
\node{G/ H} \arrow[2]{e,t}{\hat{g}_{1}} \node[2]{G/ K} \arrow[2]{e,t}{\hat{g}_{2}} \node[2]{G/L}
\end{diagram}
\]
and $\psi: G/H \times \Delta[1] \longrightarrow X$ is an equivariant $1$-simplex of type $H$  as described below. Let $x$ be a $2$-simplex in the Kan complex $X^H$ determined by the compatible pair of $1$-simplices $(g_1\phi_2^{\prime}, ~~, \phi_1^{\prime})$ so that $ \partial_0x = g_1\phi_2^{\prime}$ and $\partial_2x= \phi_1^{\prime}$. Then $\psi$ is given by $\psi (eH, \Delta_1) = \partial_1x$.
\end{Def}
It is proved in \cite{ms} that the composition is well defined.

For a version of fundamental groupoid of a $G$-space, we refer \cite{luck} and \cite{mm}.

Observe that if $X$ is a $G$-simplicial set then $S|X|$ is a $G$-Kan complex, where for any space $Y$, $SY$  denotes the total singular complex and for any simplicial set $X$, $|X|$ denotes the geometric realization of $X$.
\begin{Def}
For any $G$-simplicial set $X$, we define the fundamental groupoid $\Pi X$ of $X$ by $\Pi X := \Pi S|X|.$
\end{Def}

\begin{Rem}\lb{morphism}
If G is trivial then $\Pi X$ reduces to fundamental groupoid $\pi X$ of a simplicial set X. Again, for a fixed H, the objects $x_H$ together with the morphisms $x_{H}\rightarrow y_{H}$ with identity in the first factor, constitute a subcategory of $\Pi X$ which is precisely the fundamental groupoid $\pi X^H$ of $X^{H}$. Moreover, a morphism $[\hat{g},\phi]$ from $x_H$ to $y_K$, corresponds to the morphism $[\overline{\phi}]$ in the fundamental groupoid $\pi X^H$ of $X^H$ from $x_H^{\prime}$  to $ay_K^{\prime}$, where $\overline{\phi}$ is as in \ref{correspondence}. Suppose $\xi$ is a morphism in $\pi X^H$ from $x$ to $y$ given by a homotopy class $[\overline{\omega}],$ where $\overline{\omega}:\Delta[1]\lgra X^H$ represents the $1$-simplex in $X^H$ from $x$ to $y$. Let $x_H$ and $y_H$ be the objects in $\pi X^H$ defined respectively by $$x_H(eH,\Delta_0)=x,~~y_H(eH,\Delta_0)=y.$$ Then we have a morphism $[id, \omega]:x_H\lgra y_H$ in $\Pi X$, where $\omega(eH,\Delta_1) = \overline{\omega}(\Delta_1).$  We shall denote this morphism corresponding to $\xi$ by $b\xi.$
\end{Rem}

\begin{Def}
An equivariant local coefficients on a $G$-simplicial set $X$ is a contravariant functor from $\Pi X$ to the category $\mathcal{A}$b of abelian groups.
\end{Def}

Next, we briefly describe  the simplicial version of Bredon-Illman cohomology with local coefficients as introduced in \cite{ms}.

Let $X$ be a $G$-simplicial set and $M$ an equivariant local coefficients on $X$. For each equivariant $n$-simplex $\sigma:G/H\times \Delta[n]\rightarrow X,$ we associate an equivariant $0$-simplex $\sigma _{H}:G/H\times\Delta[0]\rightarrow X$ given by
$$\sigma_H= \sigma\circ (id\times \delta_{(1,2,\cdots,n)}),$$
where $\delta_{(1,2, \cdots ,n)}$ is the composition
$$\delta_{(1,2, \cdots ,n)}: \Delta[0]\stackrel{\delta_1}{\rightarrow} \Delta[1]\stackrel{\delta_2}{\rightarrow}\cdots \stackrel{\delta_ n}{\rightarrow} \Delta[n].$$
The $j$-th face of $\sigma$ is an equivariant $(n-1)$-simplex of type $H$, denoted by $\sigma^{(j)}$, and is defined by
$$ \sigma^{(j)}= \sigma \circ (id \times \delta_j), ~0\leq j\leq n.$$
\begin{Rem}\lb{initial}
Note that $\sigma^{(j)}_{H}=\sigma_{H}\mbox{ for }j> 0,$ and
$$\sigma^{(0)}_H = \sigma \circ (id \times \delta_{(0,2,\cdots, n)}).$$
\end{Rem}

Let $C^{n}_{G}(X;M)$ be the group of all functions $f$ defined on equivariant $n$-simplexes $\sigma:G/H\times\Delta[n]\rightarrow X$ such that $f(\sigma)\in M(\sigma _{H})$ with $f(\sigma)=0,$ if $\sigma$ is degenerate. We have a morphism $\sigma_*=[id,\alpha]$ in $\Pi X$ from $\sigma_{H}$ to $\sigma^{(0)}_{H}$ induced by $\sigma$, where $\a : G/H \times \Delta[1] \lgra X$ is given by $\a = \sigma \circ (id \times \delta_{(2,\cdots ,n)}).$ Define a homomorphism $$\delta:C_{G}^{n}(X;M)\rightarrow C_{G}^{n+1}(X;M)$$
$$f\mapsto \delta f$$
where for any equivariant $(n+1)$-simplex $\sigma $ of type $H$,
$$(\delta f)(\sigma)=M(\sigma_{*})(f(\sigma^{(0)}))+\Sigma_{j=1}^{n+1}(-1)^{j}f(\sigma^{(j)}).$$
A routine verification shows that $\delta\circ \delta=0.$ Thus $\{C_{G}^{*}(X;M),\delta \}$ is a cochain complex. We are interested in a subcomplex of this cochain complex as defined below.

Let $\eta:G/H\times \Delta[n]\rightarrow X$ and $\tau:G/ K\times \Delta[n]\rightarrow X$ be two equivariant $n$-simplexes. Suppose there exists a $G$-map $\hat{g} :G/H\lgra G/K,~ g^{-1}Hg\subseteq K,$ such that $\tau \circ (\hat{g}\times id ) = \eta.$ Then $\eta$ and $\tau $ are said to be compatible under $\hat{g}$. Observe that if $\eta$ and $\tau$ are compatible as described above then $\eta$ is degenerate if and only if $\tau$ is degenerate. Moreover notice that in this case, we have a morphism
$[\hat{g},k]:\eta_{H}\rightarrow \tau_{K}$ in $\Pi X$, where $k = \eta_H \circ (id \times \sigma_0),$ where $\sigma_0:\Delta[1]\lgra \Delta[0]$ is the map as described in Section 2. Let us denote this induced morphism by $g_*$.

\begin{Def}
We define $S_{G}^{n}(X;M)$ to be the subgroup of $C_{G}^{n}(X;M)$ consisting of all functions f such that if $\eta$ and $\tau$ are equivariant n-simplexes in X which are compatible under $\hat{g} $ then $f(\eta)=M(g_{*})(f(\tau))$.
\end{Def}

If $f\in S_{G}^{n}(X;M)$ then one can verify that $\delta f \in S_{G}^{n+1}(X;M).$

Thus we have a cochain complex $S_G(X;M) = \{S^n_G(X;M), \delta\}.$

\begin{Def}
Let X be a G-simplicial set with equivariant local coefficients M on it. Then the $n$-th Bredon-Illman cohomology of $X$ with local coefficients  $M$ is defined by $$H^n_G(X;M)=H^{n}(S_{G}(X;M)).$$
\end{Def}

\section{Equivariant Eilenberg Theorem}
Let $X$ be a one vertex Kan complex. For any $x\in X_1$, we denote by [x] the element of $\pi=\pi_1(X,v)$ represented by the $1$-simplex $x$ where $v$ is the unique vertex of $X$. Recall that (\cite{git}, \cite{gugg}) the universal covering complex $\widetilde{X}$ of X is defined as follows:
$$\widetilde{X}_n=\pi\times X_n$$
with the face maps
$$\partial_i(\gamma, x)=(\gamma, \partial_ix),~~ 0< i \leq n,~~~x\in X_n,~~\gamma\in\pi$$
$$\partial_0(\gamma, x)=([\partial_{(2,3,\cdots ,n)}x]\gamma, \partial_0x),$$ where $\partial_{(2, 3, \cdots, n)}x = \partial_2\partial_3 \cdots \partial_n x$. The degeneracy maps are
$$s_i(\gamma , x)=(\gamma, s_i x)~~~ 0\leq i\leq n.$$
Then $p:\widetilde{X}\lgra X,$ $p$ being the first projection, has the usual properties of universal covering. Any map $f:X\lgra Y$ of such complexes induces a map $\tilde{f} : \widetilde{X}\lgra \widetilde{Y}$ by $\tilde{f}(\gamma, x) = (f_*(\gamma), f(x)),$ where $f_* :\pi_1 (X)\lgra \pi_1 (Y)$ is the homomorphism of fundamental groups induced by $f$.
\begin{Rem}\lb{path}
We note that given any two $0$-simplexes $x_1=(\gamma_1, v)$ and $x_2=(\gamma_2, v)$ in $\widetilde{X},$ there is a unique homotopy class of $1$-simplexes $\omega$ such that $\partial_1\omega=x_1,~~\partial_0\omega= x_2,$ as $\widetilde{X}$ is simply connected. We may represent this class by $\omega = (\gamma_1, \omega_2\omega_1^{-1})$ where $\omega_i$ represents $\gamma_i$, $i=1,~2.$
\end{Rem}
The fundamental group $\pi_1(X)$ operates on $\widetilde{X}$ freely by
$$(\sigma, (\gamma, x))\mapsto (\gamma\sigma^{-1}, x).$$
This action is natural with respect to maps of complexes and an analogue of Eilenberg Theorem holds. The purpose of this section is to prove an equivariant version of this result. We define a contravariant functor from the category of canonical orbits to the category of one vertex Kan complexes as follows.

Let $X$ be a one vertex $G$-Kan complex. We denote the $G$-fixed vertex by $v$. For any subgroup $H$ of $G$, let $$p_H:\widetilde{X^H}\rightarrow X^H$$ be the universal cover of $X^H$. The left translation $a:X^K\rightarrow X^H,$ corresponding to a G-map $\hat{a}:G/H\rightarrow G/ K$, $a^{-1}Ha\subseteq K,$ induces a simplicial map $\tilde{a}:\widetilde{X^K}\rightarrow \widetilde{X^H}$ such that $p_H\circ   \tilde{a}=a\circ p_K$. This defines an $O_G$-Kan complex $\widetilde{X}$ by setting $\widetilde{X}(G/H)=\widetilde{X^H}$ and $\widetilde{X}(\hat{a})=\tilde{a}$. This is called the universal $O_G$-covering complex of $X$. This is simplicial analogue of $O_G$-covering space as introduced in \cite{mm}. We refer \cite{luck} for a more general version, called 'universal covering functor'. For any subgroup $H$, let $\tilde{v}^H\in \widetilde{X^H}$ denote the zero simplex $(e_H, v),$ where $e_H$ is the identity element of $\underline{\pi}X(G/H) = \pi_1(X^H, v).$ Note that the map $\tilde{a}$ induced by $a:X^K\lgra X^H$ maps $\tilde{v}^K$ to $\tilde{v}^H.$

The natural actions of $\underline{\pi}X(G/H) =\pi_1(X^H, v)$ on $\widetilde{X}(G/H)= \widetilde{X^H}$ as $H$ varies over subgroups of $G$, define an action of the $O_G$-group $\underline{\pi}X$ on $\widetilde{X}$.

Suppose $M$ is an equivariant local coefficients on $X$. We have an abelian $O_G$-group $M_0$ associated to $M$ as described below.

For any subgroup $H$ of $G$, let $v_{H}$ be the object of type $H$ in $\Pi X$ defined by
$$v_{H}:G/H\times \Delta[0]\rightarrow X,$$
$$(eH,\Delta_0)\longmapsto v.$$
Then for any morphism $\widehat{g}:G/H\rightarrow G/ K$ in $O_{G},$ given by a subconjugacy relation $g^{-1}Hg\subseteq K$, we have a morphism  $[\widehat{g},k]:v_{H}\rightarrow v_{K}$  $\Pi X$, where $k :G/H\times \Delta[1] \lgra X$ is given by $k(eH, \Delta_1)= s_0v$. Define $M_{0}:O_{G}\rightarrow Ab$ by $M_{0}(G/H)=M(v_{H})$ and $M_{0}(\widehat{g})=M[\widehat{g},k]$. The abelian  $O_{G}$-group $M_0$ comes equipped with a natural action of the $O_{G}$ group $\underline{\pi}X$ as described below.

Let $\alpha=[\phi^{\prime}]\in \underline{\pi}X(G/H)=\pi_{1}(X^{H},v)$. Then the morphism $[id,\phi]:v_{H}\rightarrow v_{H}$ where $\phi(eH,\Delta_1)= \phi^{\prime}(\Delta_1),$ is an equivalence in the category $\Pi X$. This yields a group homomorphism
$$b:\pi_{1}(X^{H},v)\rightarrow Aut_{\Pi X}(v_{H}), ~\a=[\phi^{\prime}] \mapsto b(\a)=[id,\phi].$$ We remark that the composition of the fundamental group $\pi_{1}(X^{H},v)$ coincides with the morphism composition of $\Pi X$, contrary to the usual topological composition in the fundamental group. The composition of the map $b$ with the group homomorphism $Aut_{\Pi X}(v_{H})\rightarrow Aut_{\mathcal{A} b}(M(v_{H}))$ which sends $\a\in Aut_{\Pi X}(v_{H})$ to $[M(\a)]^{-1}$ defines the action of $\pi_{1}(X^{H},v)$ on $M_{0}(G/H)$. It is routine to check that this action is natural with respect to morphism of $O_{G}$.

We define cohomology groups of $\widetilde{X}$ with coefficients in $M_0$ using invariance of actions of $\underline{\pi}X$ on $\widetilde{X}$ and $M_0$ as follows.

We have a chain complex $\left\lbrace \underline{C}_n(\widetilde{X}),\underline{\partial}_n\right\rbrace$ in the abelian category $\mathcal{C}_G$ defined by
$$\underline{C}_n(\widetilde{X})(G/H):=C_n(\widetilde{X}^H;\mathbb{Z})$$ for every object $G/H$ and for every morphism $\hat{a} :G/H\lgra G/K$, $a^{-1}Ha\subseteq K,$ $\underline{C}_n(\hat{a}):=\tilde{a}_{\#}:C_n(\widetilde{X}^K;\mathbb{Z})\rightarrow C_n(\widetilde{X}^H;\mathbb{Z})$, where $C_n(\widetilde{X}^H;\mathbb{Z})$ denotes the free abelian group generated by the non-degenerate $n$-simplexes of $\widetilde{X}^H$. The boundary map $\underline{\partial}_n:\underline{C}_n(\widetilde{X})\rightarrow \underline{C}_{n-1}(\widetilde{X})$ is the natural transformation, defined by $\underline{\partial}_n(G/H):=\partial_n$ where $\partial_n:C_n(\widetilde{X}^H;\mathbb{Z})\rightarrow C_{n-1}(\widetilde{X}^H;\mathbb{Z})$ is the ordinary boundary map. Note that $\underline{\pi}X$ acts on the chain complex $\{\underline{C}_n(\widetilde{X}), \underline{\partial}_n\}$ via its action on $\widetilde{X}$. We now form the cochain complex $\left\lbrace Hom_{\underline{\pi}X}(\underline{C}_n(\widetilde{X}),M_0),\delta^n \right\rbrace$, where $Hom_{\underline{\pi}X}(\underline{C}_n(\widetilde{X}),M_0)$ consists of all natural transformations $\underline{C}_n(\widetilde{X}) \lgra M_0$ respecting the actions of $\underline{\pi}X$ and $\delta^nf=f\circ \partial_{n+1}$. Then the $n^{th}$ equivariant cohomology of $\widetilde{X}$ with coefficients in $M_0$ is defined by
$$H^n_{\underline{\pi}X,G}(\widetilde{X};M_0):= H_n(Hom_{\underline{\pi}X}(\underline{C}_{\#}(\widetilde{X}),M_0)).$$

\begin{Thm}
Let X be a one vertex G-Kan complex and M be an equivariant local coefficients on X. Then
$$H^n_G(X ; M)\cong H^n_{\underline{\pi}X,G}(\widetilde{X};M_0).$$
\end{Thm}

\begin{proof} Recall that for any two $0$-simplexes $x,~y \in \widetilde{X^H}$ of the universal cover of the $H$-fixed point complex $X^H$, there is a unique homotopy class of $1$-simplexes $\omega $ with $\partial_1\omega = x$ and $\partial_0\omega =y$. Let us denote this class by $\widetilde{\xi}_H(x, y).$
In particular, if $x= \tilde{v}^H$, then we shall write $\widetilde{\xi}(\tilde{v}^H, y)$ simply by $\widetilde{\xi}_H(y).$ Upon projecting $\widetilde{\xi}_H(x, y)$ via $p_H$ we get an element $\xi_H(x, y) \in \pi_1(X^H, v).$  By Remark \ref{morphism}, $\xi_H(x, y)$ corresponds to an automorphism $b\xi_H(x,y)$ of $v_H$ in $\Pi X.$ As before $\xi_H(x,y)$ will be denoted by $\xi_H(y)$ when $x= \tilde{v}^H.$ Define a map
$$\phi : S^n_G(X;M) \lgra Hom_{\underline{\pi}X}(\underline{C}_n(\widetilde{X}), M_0)$$ as follows. Let $f\in S^n_G(X; M)$ and $y$ be a non-degenerate $n$-simplex in $\widetilde{X^H}.$ Let $\sigma $ be the equivariant $n$-simplex of type $H$ in $X$ such that $\sigma^{\prime} = p_H\circ \overline{y},$ where $\overline{y} :\Delta[n] \lgra \widetilde{X^H}$ is the simplicial map with $\overline{y}(\Delta_n)=y.$ Then $\phi(f) \in Hom_{\underline{\pi}X}(\underline{C}_n(\widetilde{X}), M_0)$ is given by
$$\phi (f)(G/H)(y) = M(b\xi_H(\partial_{(1,2,\cdots ,n)}y))f(\sigma).$$ Recall that $f(\sigma) \in M(\sigma_H)$ and $\sigma_H$ in this case coincides with $v_H.$

We check that $\phi (f)(G/H)$ is equivariant with respect to the respective actions of $\underline{\pi}X(G/H)$ on $\underline{C}_n(\widetilde{X})(G/H)$ and on $M_0(G/H).$
Let $u\in \underline{\pi}X (G/H), y \in \widetilde{X^H_n}$ and $\sigma$ be the equivariant $n$-simplex determined by $y$ as above. Then
$$\phi (f)(G/H)(uy) = M(b\xi_H(\partial_{(1,2,\cdots ,n)}uy))f(\tau),$$ where $\tau^{\prime} =p_H\circ \overline{uy}.$ By the definition of the action of $\underline{\pi}X (G/H)$ on $C_n(\widetilde{X^H}; \mathbb{Z})$, we have $p_H\circ \overline{uy} = p_H\circ \overline{y},$ hence $\tau^{\prime} = \sigma^{\prime}.$ It follows that
$$\phi (f)(G/H)(uy) = M(b\xi_H(\partial_{(1,2,\cdots ,n)}uy))f(\sigma).$$ Now write
$\widetilde{\xi}_H(\partial_{(1,2,\cdots ,n)}uy)$ as a composition
$$\widetilde{\xi}_H(u\tilde{v}^H,\partial_{(1,2,\cdots ,n)}uy)\circ \widetilde{\xi}_H(u\tilde{v}^H)$$ of morphisms in the fundamental groupoid of $\widetilde{X^H}$. Observe that by Remark \ref{path}, $\xi_H(u\tilde{v}^H) = u^{-1}$ and $\xi_H(u\tilde{v}^H,\partial_{(1,2,\cdots ,n)}uy)= \xi_H(\partial_{(1,2,\cdots ,n)}y).$ Therefore
$$M(b\xi_H(\partial_{(1,2,\cdots ,n)}uy)) = M(bu)^{-1}\circ M(b\xi_H(\partial_{(1,2, \cdots ,n)}y)).$$ Thus $\phi (f)(G/H)(uy) = M(bu)^{-1} \phi (f)(G/H)(y).$ It follows from the definition of the action of $\underline{\pi}X(G/H)$ on $M_0(G/H)$ that $\phi (f)(G/H)$ is equivariant.

To check that $\phi (f)(G/H) : \underline{C}_n(\widetilde{X}) \lgra M_0$ is natural, we have to show that $$M_0(\hat{g})\circ\phi (f)(G/K) = \phi (f)(G/H)\circ\tilde{g}_{\#}$$ whenever $g^{-1}Hg \subseteq K$. Recall that by definition of $M_0$, $M_0(\hat{g}) =M(v_H\xrightarrow{[\hat{g}, k]} v_K)$, where $k :G/H \times \Delta[1]\lgra X$ is given by $k(eH,\Delta_1) = s_0v.$ Let $y\in \widetilde{X^K_n}$ and $g^{-1}Hg \subseteq K.$ Let $\tau$ be an equivariant $n$-simplex of type $K$ in $X$ such that $\tau^{\prime} = p_K\circ \overline{y}$. Then
\begin{equation*}
 \begin{split}
 & M_0(\hat{g})\circ\phi (f)(G/K)(y)\\
=& M(v_{H}\xrightarrow{[\hat{g},k]} v_K)\circ M(b\xi_K(\partial_{(1,2,\cdots,n)}y))f(\tau)\\
=& M(v_{H}\xrightarrow{[\hat{g},k]} v_K)\circ M([id_{G/K},\omega])f(\tau)\\
=& M([id_{G/H}, \omega]\circ [\hat{g},k])f(\tau),
 \end{split}
\end{equation*}
where as in Remark \ref{morphism}, $\omega$ is the equivariant $1$-simplex of type $K$ in $X$ such that $\omega^{\prime}$ represents $\xi_K(\partial_{(1,2,\cdots,n)}y).$
On the other hand,
\begin{equation*}
 \begin{split}
 & \phi (f)(G/H)\circ \tilde{g}_{\#}(y)\\
=& \phi (f)(G/H)(\tilde{g} y)\\
=& M(b\xi_H(\partial_{(1,2, \cdots ,n)}\tilde{g}y))f(\sigma)
 \end{split}
\end{equation*}
where $\sigma^{\prime} = p_H\circ \overline{\tilde{g}y} = p_H\circ \tilde{g}\circ \overline{y} = g\circ p_K\circ\overline{y} =g\circ \tau^{\prime}.$ In particular, $\sigma$ and $\tau$ are compatible $n$-simplexes. Thus
$$\phi (f)(G/H)\circ \tilde{g}_{\#}(y) = M(b\xi_H(\partial_{(1,2, \cdots ,n)}\tilde{g}y))\circ M(g_*)f(\tau).$$
Note that $v$ is the only vertex in $X$ which is $G$-fixed and hence $g_*$ is a morphism from $v_H$ to $v_K$ and is given by $[\hat{g}, k]$ where $k= v_H\circ (id_{G/H} \times \sigma_0)$. Now observe that
$\xi_H(\partial_{(1,2, \cdots ,n)}\tilde{g}y)= \xi_H(\tilde{g}\partial_{(1,2, \cdots ,n)}y)$ can be represented by $g\omega^{\prime}.$ As a consequence we may write
$$b\xi_H(\partial_{(1,2, \cdots ,n)}\tilde{g}y) = [id_{G/H}, \omega\circ (\hat{g}\times id_{\Delta[1]})].$$
Therefore
\begin{equation*}
 \begin{split}
 & \phi (f)(G/H)\circ \tilde{g}_{\#}(y)\\
=& M([id_{G/H}, \omega\circ (\hat{g}\times id_{\Delta[1]})]) \circ M([\hat{g},k])f(\tau)\\
=& M([\hat{g},k]\circ [id_{G/H}, \omega\circ (\hat{g}\times id_{\Delta[1]})]) f(\tau).
 \end{split}
\end{equation*}
From the definition of composition of morphism in $\Pi X$, we have
$$[id_{G/K}, \omega]\circ [\hat{g},k] = [\hat{g},k]\circ [id_{G/H}, \omega\circ (\hat{g}\times id_{\Delta[1]})].$$
Thus $\phi (f)$ is natural.

Next we check that $\phi$ is a chain map. Let $f\in S^n_G(X;M),~ y\in \widetilde{X^H_{n+1}}.$ Let $\sigma$ denotes the equivariant $(n+1)$-simplex of type $H$ corresponding to $y$ as described before. Then
\begin{equation*}
\begin{split}
&\phi(\delta f)(G/H)(y)\\
=& M(b\xi_H(\partial_{(1,2,\cdots,n+1)}y))(\delta f)(\sigma)\\
=& M(b\xi_H(\partial_{(1,2,\cdots,n+1)}y))\left\lbrace M(\sigma_*)f(\sigma^{(0)})
   +\Sigma_{j=1}^{n+1}(-1)^{j}f(\sigma^{(j)})\right\rbrace .
\end{split}
\end{equation*}
On the other hand,
\begin{equation*}
 \begin{split}
 &  \delta\phi (f)(G/H)(y)\\
=& \Sigma_{i=0}^{(n+1)}(-1)^{i}\phi (f)(G/H)(\partial_iy)\\
=& \Sigma_{i=0}^{(n+1)}(-1)^{i}
   M(b\xi_H(\partial_{(1,2,\cdots,n)}\partial_iy))f(\sigma^{(i)})\\
=& M(b\xi_H(\partial_{(0,2,\cdots,n+1)}y))f(\sigma^{(0)})
   +\Sigma_{i=1}^{n+1}(-1)^{i}M(b\xi_H(\partial_{(1,2,\cdots,n+1)}y))f(\sigma^{(i)}).
\end{split}
\end{equation*}
Now note that since $\widetilde{X^H}$ is simply connected the morphism $\xi_H(\partial_{(0,2,\cdots,n+1)}y)$ in $\pi X^H$ can be factored as
$$\xi_H(\partial_{(1,2,\cdots,n+1)}y, \partial_{(0,2,\cdots,n+1)}y)\circ \xi_H(\partial_{(1,2,\cdots,n+1)}y)$$
and $b\xi_H(\partial_{(1,2,\cdots,n+1)}y, \partial_{(0,2,\cdots,n+1)}y)$ is precisely
the morphism $\sigma_*$. Therefore
$$b\xi_H(\partial_{(0,2,\cdots,n+1)}y) = \sigma_*\circ b\xi_H(\partial_{(1,2,\cdots,n+1)}y).$$
Hence $\phi(\delta f) = \delta\phi (f).$

To show that $\phi$ is a chain isomorphism define a map  $$\psi:Hom_{\underline{\pi}X}(\underline{C}_n(\tilde{X}),M_0) \rightarrow  C_G^n(X;M)$$ as follows. Let $f\in Hom_{\underline{\pi}X}(\underline{C}_n(\tilde{X}),M_0)$ and $\sigma$ be a non-degenerate equivariant $n$-simplex in X of type $H.$ Choose an $n$-simplex $y$ in $\widetilde{X^H}$ such that $p_H(y) = \sigma (eH, \Delta_n).$ Then $\psi (f)$ is given by
$$\psi (f)(\sigma ) = M(b\xi_H(\partial_{(1,2, \cdots, n)}y))^{-1}f(G/H)(y).$$

Suppose $z$ is another $n$-simplex in $\widetilde{X^H}$ such that $p_H(z)= \sigma (eH, \Delta_n).$ Since $\pi_1(X^H,v)$ acts transitively on each fibre of $p_H:\widetilde{X^H}\lgra X^H$, there exists an element $u\in \pi_{1}(X^H,v)$ such that $u y= z$ and hence $u\partial_{(1,2, \cdots, n)}y = \partial_{(1,2, \cdots, n)}z.$ Thus
\begin{equation*}
 \begin{split}
 & M(b\xi_H(\partial_{(1,2, \cdots, n)}z))^{-1}f(G/H)(z)\\
=& M(b\xi_H(\partial_{(1,2, \cdots, n)}uy))^{-1}f(G/H)(u y)\\
=& M(b\xi_H(\partial_{(1,2, \cdots, n)}y))^{-1}M(bu)^{-1}f(G/H)(y)\\
=& M(b\xi_H(\partial_{(1,2, \cdots, n)}y))^{-1}f(G/H)(y).
 \end{split}
\end{equation*}
The last equality follows from the observation
$$M(b\xi_H(\partial_{(1,2,\cdots ,n)}uy)) = M(bu)^{-1}\circ M(b\xi_H(\partial_{(1,2, \cdots ,n)}y)),$$ while showing $\phi$ takes any cocycle in $S^n_G(X;M)$ into $Hom_{\underline{\pi}X}(\underline{C}_n(\tilde{X}),M_0),$ in the first part of the proof. Thus the map $\psi$ is well defined.

We claim that $\psi (f) \in S^n_G(X;M)$ for any $f\in Hom_{\underline{\pi}X}(\underline{C}_n(\tilde{X}),M_0).$ Let $a^{-1}Ha\subseteq K$ and $\sigma :G/H\times \Delta[n] \lgra X$ and $\eta: G/K\times \Delta[n]\lgra X$ be equivariant $n$-simplexes such that $\eta \circ (\hat{a}\times id) = \sigma$ so that they are compatible. We need to show that $\psi (f)(\sigma)= M(a_*)\psi (f)(\eta).$
Let $y\in \widetilde{X^K}$ be such that $p_K(y) = \eta (eK, \Delta_n).$ Then the $n$-simplex $\tilde{a}y \in \widetilde{X^H_n}$ is such that $$p_H(\tilde{a}y)=ap_K(y)=a\eta(eK, \Delta_n) = \eta (aK, \Delta_n)= \sigma (eH,\Delta_n).$$
By our choice, we have
$$\psi (f)(\eta) = M(b\xi_K(\partial_{(1,2, \cdots, n)}y))^{-1}f(G/K)(y)$$
and $$\psi (f)(\sigma) = M(b\xi_H(\partial_{(1,2, \cdots, n)}\tilde{a}y))^{-1}f(G/H)(\tilde{a}y).$$
Since $f:\underline{C}_n(\tilde{X})\lgra M_0$ is natural, we have
$$f(G/H)(\tilde{a}y) =M_0(\hat{a})(f(G/K)(y)).$$
In the first part of the proof we have observed that
$$a_*\circ b\xi_H(\partial_{(1,2, \cdots ,n)}\tilde{a}y) =b\xi_K(\partial_{(1,2, \cdots, n)}y)\circ a_*.$$ Moreover, recall that $M_0(\hat{a}) = M(a_*).$ Therefore
\begin{equation*}
 \begin{split}
 & M(a_*)\psi (f)(\eta)\\
=& M(a_*)M(b\xi_K(\partial_{(1,2, \cdots, n)}y))^{-1}f(G/K)(y)\\
=& M(b\xi_K(\partial_{(1,2, \cdots, n)}y)^{-1}\circ a_*)f(G/K)(y)\\
=& M(a_*\circ b\xi_H(\partial_{(1,2, \cdots, n)}\tilde{a}y))f(G/K)(y)\\
=& M(b\xi_H(\partial_{(1,2, \cdots, n)}\tilde{a}y))^{-1}M(a_*)f(G/K)(y)\\
=& M(b\xi_H(\partial_{(1,2, \cdots, n)}\tilde{a}y))^{-1}M_0(\hat{a})f(G/K)(y)\\
=& \psi (f)(\sigma).
 \end{split}
\end{equation*}
It is routine to check that $\psi$ is the inverse of $\phi$. This completes the proof of the the theorem.
\end{proof}
\section{Equivariant Serre spectral sequence} The aim of this last section is to derive a version of Serre spectral sequence. To do this we give an alternative description of equivariant simplicial cohomology with local coefiicients in terms of cohomology of small categories.

Let $G$ be a discrete group and $X$ a $G$-Kan complex. Then we have a category $\Delta_G(X)$ described as follows. Its objects are $G$-simplicial maps $\sigma:G/H\times \Delta[n]\rightarrow X$ and a morphism from $\sigma:G/ H\times \Delta[n]\rightarrow X$  to $\tau:G/ K\times \Delta[m]\rightarrow X$ is a pair $(\hat{g},\alpha)$ where $\hat{g}:G/ H\rightarrow G/ K$ is a $G$-map and $\alpha:\Delta[n]\rightarrow \Delta[m]$ is a simplicial map such that $\tau\circ (\hat{g}\times \alpha)=\sigma.$ There is a canonical functor $v_X:\Delta_G(X)\rightarrow \Pi(X)$ which sends $\sigma:G/H\times\Delta[n]\rightarrow X$ to $\sigma_H=\sigma\circ (id\times \delta_{(1,2,\cdots,n)}).$ For a morphism $(\hat{g},\alpha)$ in $\Delta_G(X)$, $v_X(\hat{g},\alpha):\sigma_H\rightarrow \tau_K$ is the morphism $[\hat{g},\phi]$ in $\Pi (X)$ where $\phi:G/H\times \Delta[1]\rightarrow X$ is an equivariant $1$-simplex of type $H$ obtained as follows. Suppose $\tau\circ(id\times \delta_{(1,\cdots,\widehat{\alpha(0)}\cdots,m)})(eK,\Delta_1)=\omega\in X^K.$  Let $x$ be a $2$-simplex in $X^K$ determined by the compatible pair of $1$-simplices $(~~, s_1\partial_1\omega, \omega).$ Then $\phi$ is given by $\phi(eH,\Delta_1)=g(\partial_0 x).$

If $X$ is any $G$-simplicial set then we define $\Delta_G(X) = \Delta_G(S|X|).$

For a small category $\mathcal{C},$ let $\mathcal{A}b(\mathcal{C})$ be the category of all contravariant functors from $\mathcal{C}$ to $\mathcal{A}b$ with morphisms natural transformations of functors.
\begin{Def}
A functor $M\in \mathcal{A}b(\Delta_G(X))$ is said to be $G$-local if
$$M=v_X^*M^{\prime}=M^{\prime}\circ v_X$$ for some $M^{\prime}\in \mathcal{A}b(\Pi (X)).$ For a $G$-local coefficients  $M$, the equivariant cohomology of $X$ with coefficients $M$ is defined to be
$$h_G^*(X;M):= H^*(\Delta_G(X);M),$$
where the right hand side denotes the cohomology of the category $\Delta_G(X)$, in the sense of \cite{qui}.
\end{Def}
\begin{Thm}\lb{iso} Let $X$ be a $G$-simplicial set and $M$ be an equivariant local coefficients on $X$. Then there is an isomorphism $$H^*_G(X;M)\cong h_G^*(X;M). $$
(On the right we identify $M$ with $v_X^\ast(M)$).
\end{Thm}

\begin{proof}
Let $N\mathcal{C}$ denote the nerve of a small category $\mathcal{C}$. Then as in (\cite{moer}) we let $\tilde{X}$ be the bisimplicial set whose $(p,q)$ simplices are triples $(u,\alpha,\sigma)$ where
$$\begin{array}{l}
u=([n_0]\xrightarrow{u_1}[n_1]\rightarrow\cdots\xrightarrow{u_p}[n_p])\in N_p(\Delta) \\
\alpha=(G/H_0\xrightarrow{\alpha_1}G/H_1\rightarrow\cdots\xrightarrow{\alpha_q}G/H_q)\in N_q(O_G)\\
\sigma:G/H_q\times\Delta[n_p]\rightarrow X\mbox{ is a $G$-simplicial map.}
\end{array}$$

The face and degeneracy maps on $\tilde{X}$ are induced from those on $N(\Delta)$ and $ N(O_G)$. Then $$\mbox{diagonal}(\tilde{X})\cong N(\Delta_G(X)).$$ To every $(u,\alpha,\sigma)\in\tilde{X}^{p,q}$ associate a $G$-simplicial map, $$\overline{\sigma}=\sigma\circ(\alpha_q\circ\cdots\circ\alpha_1\times u_p\circ\cdots\circ u_1):G/{H_0}\times\Delta[n_0]\rightarrow X.$$

Let $C^{p,q}(X;M)$ denote the set of all functions on $\tilde{X}^{p,q}$ which sends an element $(u,\alpha,\sigma)$ of $\tilde{X}^{p,q}$ to an element of $M(v_X(\overline{\sigma}))$. It follows quite easily that $C^{p,q}(X;M)$ is a bicomplex with obvious differentials $d_h$ and $d_v$ induced from the face maps of $\tilde{X}$. Denote the total complex of $C^{\bullet\bullet}(X,M)$ by $\mbox{Tot}\,C^{\bullet\bullet}(X;M)$.
Let $\mbox{diag}\,C^{\bullet\bullet}(X;M)$ be the cochain complex whose $p^{th}$ group is $C^{p,p}(X;M)$ and differential is $d_hd_v$. Then by a result of Dold and Puppe (\cite{dp}) we have $$H^n(\mbox{Tot}\,C^{\bullet\bullet}(X;M))\cong H^n(\mbox{diag}\,C^{\bullet\bullet}(X;M)).$$ Now $C^{p,p}(X;M)$ can be interpreted as the set of all functions on
$N(\Delta_G(X))$ which sends a $p$-simplex $(\tau_0\rightarrow\tau_1\rightarrow\cdots\rightarrow \tau_p)$ to an element of $M(v_X(\tau_0))$ and the differential on $C^{p,p}(X;M)$ is just the differential induced from the face maps of $N_p(\Delta_G(X))$. Hence, $$ H^n(\mbox{diag}\,C^{\bullet\bullet}(X;M))\cong
H^n(\Delta_G(X);v_X^\ast M)=h^n_G(X;M).$$

Recall that the spectral sequence associated to the $p$-filtration of the bicomplex $C^{\bullet\bullet}(X;M)$ converges to the cohomology of the total complex. Now proceeding as in \cite{mp}, we may compute the $E_1$ term of the spectral sequence. It turns out that

$$\begin{array}{lcll} E_1^{p,q}= H^q(C^{p,\bullet}(X;M)) & \cong & \Pi_{u\in N_p(\Delta)}S^{n(u)}_G(X;M) & \mbox{if }q=0  \\
& = & 0 & \mbox{if }q>0,
\end{array}$$
where $S^{n(u)}_G(X;M)$ is a copy of $S^{n_p}_G(X;M)$ for every $u = ([n_0]\ra \cdots\ra [n_p]).$ (See \cite{mp} for details.)
Thus,
$$\begin{array}{lcl}
H^p(\mbox{Tot }C^{\bullet\bullet}(X;M)) & \cong &
H^p(\Pi_{u\in N(\Delta)}S^{n(u)}_G(X;M))
\\
  &  &  \\
 & \cong & H^p(\Delta^{op},S^\bullet_G(X;M))
 \end{array}$$
where $S^\bullet_G(X;M)$ is the cosimplicial group which takes $[n]$ to $S^n_G(X;M)$ with obvious face and degeneracy maps induced from those on $\Delta$. Then we know that (\cite{moer}), $$H^p(\Delta^{op};S^\bullet_G(X;M))\cong H^p(S^\bullet_G(X;M)).$$ Hence, $$H^p(\mbox{Tot }C^{\bullet\bullet}(X;M)) \cong H^p_G(X;M).$$
\end{proof}

We are now in a position to derive the required spectral sequence.

Let $X,Y$ be $G$-simplicial sets and $f:Y\rightarrow X$ be a $G$-Kan fibration and $M$ a $G$-local coefficients  on $Y$. For $q\geq 0$ we have a contravariant functor
$$h_G^q(f,M):\Delta_G(X)\rightarrow \mathcal{A}b$$
as follows. For an object $\sigma:G/H\times \Delta[n]\rightarrow X$ of $\Delta_G(X)$, let $\sigma^*(Y)$ be the $G$-simplicial set obtained by pulling back $f$ along $\sigma$ and define $$h_G^q(f,M)(\sigma):=h_G^q(\sigma^*(Y);\tilde{\sigma}^*M),$$ where $\tilde{\sigma}:\sigma^*(Y)\rightarrow Y$ is the canonical map. We claim that $h_G^q(f,M)$ factors through $v_X$ yielding a $G$-local coefficients on $X$. To see this, first note that the following result holds (\cite{moer}, the proof of the Theorem 2.3).
\begin{Thm}
Let $f:Y\lgra X$ be a weak equivalence in $G\mathcal{S}$. Then for any $G$-local coefficients $M$ on $X$, $f$ induces an isomorphism
$$h^*_G(X; M) \cong h^*_G(Y; f^*M).$$
\end{Thm}
Recall that $f:Y\lgra X$ is a weak equivalence in $G\mathcal{S}$ if and only if $f^H: Y^H\lgra X^H$ is a weak equivalence in $\mathcal{S}$.
Next note that if $u: \Delta[m] \lgra \Delta[n]$ is a simplicial map then the map $\sigma(id\times u))^*(Y) \lgra \sigma^*(Y)$ covering $id \times u :G/K\times \Delta [m] \lgra  G/K\times \Delta [n]$ is a weak equivalence in $G\mathcal{S}$, as $f$ is a $G$-Kan fibration. The claim now follows from the above result.
\begin{Thm}
For any $G$-Kan fibration $f:Y\rightarrow X$ and a $G$-local coefficients $M$ on $Y$, there is a natural spectral sequence with $E_2$-term $E_2^{p,q}=H^p_G(X;h^p_G(f,M))$  converging to $H^{p+q}_G(Y;M).$
\end{Thm}

\begin{proof}
 The proof is parallel to the proof of Theorem 3.2, \cite{moer}. We only mention the essential steps. The $G$-Kan fibration $f:Y\rightarrow X$ induces a functor $\Delta_G(f):\Delta_G(Y)\rightarrow \Delta_G(X)$ and we have a Grothendieck spectral sequence \cite{segal}
$$H^p(\Delta_G(X);h^q(\Delta_G(f)/-;M))\Rightarrow H^{p+q}(\Delta_G(Y);M).$$
It is enough to show that the two contravariant functors $h^q(\Delta_G(Y)/-;M)$ and $h^q_G(f,M)$ from $\Delta_G(X)$ to $\mathcal{A}$b  are equivalent.
For an object $\sigma$ of $\Delta_G(X),$  $\Delta_G(f)/\sigma$ is the comma category. Objects of $\Delta_G(f)/\sigma$ are pairs $(\tau,u)$ where $\tau\in Ob(\Delta_G(Y))$ and $u:\Delta_G(f)(\tau)\rightarrow \sigma$ is a map in $\Delta_G(X).$ Morphisms from $(\tau,u)$ to $(\tau^{\prime},u^{\prime})$ are maps $\alpha:\tau\rightarrow \tau^{\prime}$ such that $u^{\prime}\Delta_G(f)(\alpha)=u.$
A direct computation shows that there is a canonical equivalence of the categories $$\Delta_G(f)/\sigma\cong\Delta_G(\sigma^*(Y)),$$ which is natural in $\sigma.$
Hence we have natural isomorphism of functors $$h^q(\Delta_G(\sigma^*(Y);\tilde{\sigma}^*M))\cong h^q(\Delta_G(f)/\sigma;M).$$ The result now follows from Theorem \ref{iso}.

\end{proof}

{\bf Goutam Mukherjee}\\
Indian Statistical Institute, Kolkata-700108, India.\\
e-mail: goutam@isical.ac.in

{\bf Debasis Sen}\\
Indian Statistical Institute, Kolkata-700108, India.\\
e-mail: dsen\_r@isical.ac.in

\end{document}